\documentclass[12pt]{article}         
\usepackage{authblk,graphicx,amsmath,amsfonts,amssymb,amscd,amsmath,breqn}  
\usepackage{amsthm}

\title{AN INTEGRAL'S JOURNEY OVER THE REAL LINE}

\author[1,2]{Robert Reynolds}
\author[1,2]{Allan Stauffer}
\affil[1]{Faculty of Science, York University}
\affil[2]{Department of Mathematics and Statistics, York University}
\affil[ ]{milver@my.yorku.ca, stauffer@yorku.ca}

\begin{document}

\maketitle

\begin{abstract}

In 1826 Cauchy presented an Integral over the real line. Al and I thought a derivation would be mighty fine. So we packed our contour integral bags that day, and we now present an analytic continuation this time.

\end{abstract}

\section{Introduction}

While reading the works of Cauchy \cite{cauchy} and Bierens de haan \cite{bdh}, we encountered a vast number of definite integrals. We found interest, in one particular definite integral which had a logarithmic and rational function integrand over the real line. Upon further reading we realized a formal derivation was not supplied in the literature. We decided to produce a formal derivation for this integral and during this process we were able to derive a generalized integral form which yielded the Cauchy integral as a special case. Another result of our work was the derivation of integrals featuring $\log(\log(x))$ integrands which were published by Malmsten \cite{malm1}. The general integral form we derived is given by

\begin{dmath}
\int_{-\infty}^{\infty}\frac{\log ^k(y) \log (1-y z)}{a^2+y^2}dy
\end{dmath}

is derived in terms of the logarithmic function, where the parameters $a,k$ and $z$ are general complex numbers subject to the restrictions given.  The closed form solution of the derived generalized integral is the product of the logarithmic function, this is unique as the the integrand too is a product of the logarithmic function. We thought these characteristics along with taking the integral over the real line strengthens the case made to study such integrals.



  The derivations follow the method used by us in \cite{reyn2}. The generalized Cauchy's integral formula is given by

\begin{equation}\label{intro:cauchy}
\frac{y^k}{k!}=\frac{1}{2\pi i}\int_{C}\frac{e^{wy}}{w^{k+1}}dw.
\end{equation}

 This method involves using a form of equation (\ref{intro:cauchy}) then multiply both sides by a function, then take a definite integral of both sides. This yields a definite integral in terms of a contour integral. A second contour integral is derived by multiplying equation (\ref{intro:cauchy}) by a function and performing some substitutions so that the contour integrals are the same.
 

\section{Derivation of the first contour integral}

We use the method in \cite{reyn2}. Here the contour is similar to Figure 2 in \cite{reyn2}. Using a generalization of Cauchy's integral formula we first replace $y$ by $\log(y)$ followed by multiplying both sides by

\begin{equation*}
\frac{\log (1-zy)}{y^2+a^2}
\end{equation*}

 to get

\begin{dmath}\label{eq1}
\frac{\log ^k(y) \log (1-y z)}{a^2+y^2}=\frac{1}{2\pi i}\int_{C}\frac{w^{-k-1} y^w \log (1-y z)}{a^2+y^2}dw
\end{dmath}

where the logarithmic function is defined in equation (4.1.2) in \cite{as}. We then take the definite integral over $y \in (-\infty,\infty)$ of both sides to get

\begin{dmath}\label{eq2}
\int_{-\infty}^{\infty}\frac{\log ^k(y) \log (1-y z)}{a^2+y^2}dx
=\frac{1}{2\pi i}\int_{-\infty}^{\infty}\int_{C}\frac{w^{-k-1} y^w \log (1-y z)}{a^2+y^2}dwdx
=\frac{1}{2\pi i}\int_{C}\int_{-\infty}^{\infty}\frac{w^{-k-1} y^w \log (1-y z)}{a^2+y^2}dxdw
=\frac{1}{2\pi i}\int_{C}\frac{1}{2} \pi  e^{\frac{i \pi  w}{2}} \left(\frac{1}{a}\right)^{1-w} w^{-k-1} \left(\log \left(a^2 z^2+1\right)-2 i
   \tan ^{-1}(a z)\right)dw
\end{dmath}

from equation (145)(27) in \cite{bdh} and the integral is valid for $a$, $z$ and $k$ complex and $-1<Re(w)<0$.

\section{Derivation of the second contour integral}

In this section we will again use the generalized Cauchy's integral formula equation (\ref{intro:cauchy}) to derive equivalent contour integrals. First we replace $y$ by $\log (a)+\frac{i \pi }{2}$ and multiply both sides by $-\frac{i \pi  \tan ^{-1}(a z)}{a}$ to get

\begin{dmath}\label{eq3}
-\frac{i \pi  \left(\log (a)+\frac{i \pi }{2}\right)^k \tan ^{-1}(a z)}{a k!}=-\frac{1}{2\pi i}\int_{C}i \pi  e^{\frac{i \pi  w}{2}} a^{w-1} w^{-k-1} \tan ^{-1}(a z)dw
\end{dmath}

where $-1<Re(w)<0$, and then using equation (\ref{intro:cauchy}) we replace $y$ by $\log (a)+\frac{i \pi }{2}$ and multiply both sides by $\frac{\pi  \log \left(a^2 z^2+1\right)}{2 a}$ to get

\begin{dmath}\label{eq4}
\frac{\pi  \left(\log (a)+\frac{i \pi }{2}\right)^k \log \left(a^2 z^2+1\right)}{2 a k!}=\frac{1}{2\pi i}\int_{C}\frac{1}{2} \pi  e^{\frac{i \pi  w}{2}} a^{w-1} w^{-k-1} \log \left(a^2
   z^2+1\right)dw
\end{dmath}

where $-1<Re(w)<0$.

\section{Derivation of the definite integral in terms of the logarithmic function}

Since the right-hand side of equation (\ref{eq2}) is equal to the sum of the right-hand sides of equations (\ref{eq3})  and (\ref{eq4}) we can equate the left-hand sides simplifying to get 

\begin{dmath}\label{eq5}
\int_{-\infty}^{\infty}\frac{\log ^k(y) \log (1-y z)}{a^2+y^2}dy=\frac{\pi  \log ^k(i a) \log (1-i a z)}{a}
\end{dmath}

where $Im(z)>0$.

\section{Derivation of a logarithmic form using substitution}

Using equation (\ref{eq5}) we replace $z$ by $p$ to form a second equation then adding these two equations to get

\begin{dmath}
\int_{-\infty}^{\infty}\frac{\log ^k(y) (\log (1-p y)+\log (1-y z))}{a^2+y^2}dy=\frac{\pi  \log ^k(i a) (\log (1-i a p)+\log (1-i a z))}{a}
\end{dmath}

Next we replace $z$ by $r \cos (t)+i r \sin (t)$ and $p$ by $-r \cos (t)+i r \sin (t)$ to get

\begin{dmath}
\int_{-\infty}^{\infty}\frac{\log ^k(y) \left(\log \left(1+r e^{-i t} y\right)+\log \left(1-r e^{i t} y\right)\right)}{a^2+y^2}dy=\frac{\pi  \log ^k(i a) (\log (1+a r (\sin (t)-i \cos
   (t)))+\log (1+a r (\sin (t)+i \cos (t))))}{a}
\end{dmath}

We then simplify by combining the logarithmic functions to get

\begin{dmath}
\int_{-\infty}^{\infty}\frac{\log ^k(y) \log \left(-r^2 y^2-2 i r y \sin (t)+1\right)}{a^2+y^2}dy=\frac{\pi  \log ^k(i a) \log \left(a^2 r^2+2 a r \sin (t)+1\right)}{a}
\end{dmath}

Since we are integrating over $y \in (-\infty,\infty)$ we can replace $y$ by $-y$ in one of the integrals in order to eliminate the imaginary part, say in the second integral to get

\begin{dmath}
\int_{-\infty}^{\infty}\frac{1}{a^2+y^2}\left(\log ^k(-y) \log (i r y \sin (t)-r y \cos (t)+1)+\log ^k(y) \log (-i r y \sin (t)-r y \cos (t)+1)\right)dy\\
=\frac{\pi  \log ^k(i a) \log \left(a^2 r^2+2 a r
   \sin (t)+1\right)}{a}
\end{dmath}

We then simplifying the logarithmic functions to get

\begin{dmath}\label{dl:eq1}
\int_{-\infty}^{\infty}\frac{1}{a^2+y^2}\left(\log ^k(-y) \log (i r y \sin (t)-r y \cos (t)+1)+\log ^k(y) \log (-i r y \sin (t)-r y \cos (t)+1)\right)dy\\
=\frac{\pi  \log ^k(i a) \log \left(a^2 r^2+2 a r
   \sin (t)+1\right)}{a}
\end{dmath}

\section{Derivation of entry 2.6.14.19 in \cite{prud}}

Using equation (\ref{dl:eq1}) and replacing $r$ by $1/r$ and setting $k=0$ simplifying we get

\begin{dmath}
\int_{-\infty}^{\infty}\frac{\log \left(1-\frac{e^{-i t} y}{r}\right)+\log \left(1-\frac{e^{i t} y}{r}\right)}{a^2+y^2}dy=\frac{\pi  \log \left(\frac{a^2+2 a r \sin
   (t)+r^2}{r^2}\right)}{a}
\end{dmath}

where $\int_{-\infty}^{\infty}\frac{\log \left(r^2\right)}{a^2+y^2}dx=\frac{\pi  \log \left(r^2\right)}{a}$ from equation (3.241.2) in \cite{grad}.
\\\\
 We then simplify the logarithmic functions and simplify to get

\begin{dmath}
\int_{-\infty}^{\infty}\frac{\log \left(r^2-2 r y \cos (t)+y^2\right)}{a^2+y^2}dy=\frac{\pi  \left(2 \log \left(a^2+2 a r \sin (t)+r^2\right)-4 \log (r)\right)}{2 a}+\frac{2 \pi  \log
   (r)}{a}
\end{dmath}

Next we replace $a$ by $z$, $r$ by $a$, $y$ by $x$, $\cos (t)$ by $b$ and $\sin (t)$ by $\sqrt{1-b^2}$ simplifying to get

\begin{dmath}\label{s1:eq1}
\int_{-\infty}^{\infty}\frac{\log \left(a^2-2 a b x+x^2\right)}{x^2+z^2}dx=\frac{\pi  \log \left(a^2+2 a \sqrt{1-b^2} z+z^2\right)}{z}
\end{dmath}

where $Im(z)>0$.

\section{Derivation of entry 4.296.2 in \cite{grad}}

Using equation (\ref{s1:eq1}) and setting $a=1$ simplifying we get

\begin{dmath}
\int_{-\infty}^{\infty}\frac{\log \left(r^2-2 r y \cos (t)+y^2\right)}{y^2+1}dy=\pi  \log \left(r^2+2 r \sin (t)+1\right)
\end{dmath}

%
%
%
%
%
%
%
%
%
%
Upon inspection of the closed form solution on the right-hand side we are able to write down the conditional solution for the definite integral given by
\begin{dmath}
    \int_{-\infty}^{\infty}\frac{\log \left(r^2-2 a y \cos (t)+y^2\right)}{y^2+1}dy= 
\begin{cases}
   \pi  \log \left(r^2+2 r \sin (t)+1\right),& \text{if } r\in\mathbb{R}, t\in\mathbb{C}\\
    \pi  \log \left(r^2+2 r |\sin (t)|+1\right),& \text{if $r,t\in\mathbb{R}$}
\end{cases}
\end{dmath}

which represents the analytic continuation of the original integral.

\section{Derivation of special cases}

Using equation (\ref{eq5}) and forming a second equation by replacing $z$ by $p+iq$. Next using this new equation we replace $p$ by $-p$ to get a third equation and add this to the second equation and taking the first partial derivative with respect to $k$ to get

\begin{dmath}\label{sec8:eq1}
\int_{-\infty}^{\infty}\frac{\log (\log (y)) \log ^k(y) \log \left(1-y \left(y \left(p^2+q^2\right)+2 i q\right)\right)}{a^2+y^2}dy=\frac{\pi  \log (\log (i a)) \log ^k(i
   a) \log \left(a^2 \left(p^2+q^2\right)+2 a q+1\right)}{a}
\end{dmath}

\subsection{Example 1}

Using equation (\ref{sec8:eq1}) and setting $k=0,a=1,p=0,q=1$ simplifying to get

\begin{dmath}
\int_{-\infty}^{\infty}\frac{\log \left(-(y+i)^2\right) \log (\log (y))}{y^2+1}dy=\pi  \log (4) \log \left(\frac{i \pi }{2}\right)
\end{dmath}

\subsection{Example 2}

Using equation (\ref{sec8:eq1}) and setting $k=1/2,a=1,p=0,q=1$ simplifying to get

\begin{dmath}
\int_{-\infty}^{\infty}\frac{\sqrt{\log (y)} \log \left(-(y+i)^2\right) \log (\log (y))}{y^2+1}dy=(1+i) \pi ^{3/2} \log (2) \left(\log \left(\frac{\pi }{2}\right)+\frac{i
   \pi }{2}\right)
\end{dmath}

\subsection{Example 3}

Using equation (\ref{sec8:eq1}) and setting $k=-1/2,a=1,p=0,q=1$ simplifying to get

\begin{dmath}
\int_{-\infty}^{\infty}\frac{\log \left(-(y+i)^2\right) \log (\log (y))}{\left(y^2+1\right) \sqrt{\log (y)}}dy=(1-i) \sqrt{\pi } \log (4) \log \left(\frac{i \pi
   }{2}\right)
\end{dmath}

\subsection{Example 4}

Using equation (\ref{sec8:eq1}) and setting $k=1,a=1,p=-1,q=1$ simplifying to get

\begin{dmath}
\int_{-\infty}^{\infty}\frac{\log (y) \log (1-2 y (y+i)) \log (\log (y))}{y^2+1}dy=\frac{1}{2} i \pi ^2 \log (5) \log \left(\frac{i \pi }{2}\right)
\end{dmath}

\subsection{Example 5}

Using equation (\ref{eq5}) taking the first partial derivative with respect to $z$ and $k$ and setting $k=0,a=2$ simplifying to get

\begin{dmath}
\int_{-\infty}^{\infty}\frac{y \log (\log (y))}{(1-i y) \left(y^2+4\right)}dy=\frac{1}{3} i \pi  \log \left(\log (2)+\frac{i \pi }{2}\right)
\end{dmath}

\subsection{Example 6}

Using equation (\ref{eq5}) taking the first partial derivative with respect to $z$ and $k$ and setting $k=-1/2,a=1,z=i$ simplifying to get

\begin{dmath}
\int_{-\infty}^{\infty}\frac{y \log (\log (y))}{(1-i y) \left(y^2+1\right) \sqrt{\log (y)}}dy=\sqrt[4]{-1} \sqrt{\frac{\pi }{2}} \log \left(\frac{i \pi }{2}\right)
\end{dmath}

\section{Summary table of integrals}

\renewcommand{\arraystretch}{2.0}
\begin{tabular}{ l  c }
  \hline			
  $f(y)$ & $\int_{-\infty}^{\infty}f(y)dy$ \\ \hline
  $\frac{\log \left(a^2-2 a b y+y^2\right)}{y^2+z^2}$ & $\frac{\pi  \log \left(a^2+2 a \sqrt{1-b^2} z+z^2\right)}{z}$  \\
  $\frac{\log \left(r^2-2 r y \cos (t)+y^2\right)}{y^2+1}$ & $\begin{cases}
   \pi  \log \left(r^2+2 r \sin (t)+1\right)& \\
    \pi  \log \left(r^2+2 r |\sin (t)|+1\right)& 
\end{cases}$  \\
  $\frac{\log ^k(y) \log (1-y z)}{a^2+y^2}$ & $\frac{\pi  \log ^k(i a) \log (1-i a z)}{a}$  \\
  $\frac{\log \left(-(y+i)^2\right) \log (\log (y))}{y^2+1}$ & $\pi  \log (4) \log \left(\frac{i \pi }{2}\right)$  \\
  $\frac{\sqrt{\log (y)} \log \left(-(y+i)^2\right) \log (\log (y))}{y^2+1}$ & $(1+i) \pi ^{3/2} \log (2) \left(\log \left(\frac{\pi }{2}\right)+\frac{i
   \pi }{2}\right)$  \\
  $\frac{\log \left(-(y+i)^2\right) \log (\log (y))}{\left(y^2+1\right) \sqrt{\log (y)}}$ & $(1-i) \sqrt{\pi } \log (4) \log \left(\frac{i \pi
   }{2}\right)$  \\
  $\frac{\log (y) \log (1-2 y (y+i)) \log (\log (y))}{y^2+1}$ & $\frac{1}{2} i \pi ^2 \log (5) \log \left(\frac{i \pi }{2}\right)$  \\
  $\frac{y \log (\log (y))}{(1-i y) \left(y^2+4\right)}$ & $\frac{1}{3} i \pi  \log \left(\log (2)+\frac{i \pi }{2}\right)$  \\
  $\frac{y \log (\log (y))}{(1-i y) \left(y^2+1\right) \sqrt{\log (y)}}$ & $\sqrt[4]{-1} \sqrt{\frac{\pi }{2}} \log \left(\frac{i \pi }{2}\right)$ 
   \\[0.3cm]
  \hline  
\end{tabular}

\section{Discussion}

In this work we looked at deriving the product of the logarithmic and rational functions over the real line. One of the interesting properties of this integral is its closed form being similar to the integrand. We formally derived a few integrals by famous mathematicians and supplied a few new results as well. We were also able to produce closed form solutions for integrals with singularities, where these results could be added to existing textbooks.

The results presented were numerically verified for both real and imaginary values of the parameters in the integrals using Mathematica by Wolfram. We considered various ranges of these parameters for real, integer, negative and positive values. We compared the evaluation of the definite integral to the evaluated Special function and ensured agreement.

\section{Conclusion}
In this paper we used our method to evaluate definite integrals using the Lerch function. The contour we used was specific to  solving integral representations in terms of the Lerch function. We expect that other contours and integrals can be derived using this method.

\section{Acknowledgments}
This paper is fully supported by the Natural Sciences and Engineering Research Council (NSERC) Grant No. 504070.


\begin{thebibliography}{999}

\bibitem{cauchy} A.L. Cauchy
\emph{Analyse transcendante. Recherche d’une formule générale quifournit la valeur de la plupart des intégrales définies connueset celle d’un grand nombre d’autres}, Annales de Mathématiques pures et appliquées, tome  17 (1826-1827), p. 84-127

\bibitem{grad} Gradshteyn I.S \& Ryzhik I.M,
\emph{Tables of Integrals, Series and Products, 6 Ed}, Academic Press (2000), USA.

\bibitem{wang} Wang, Z.X., Guo, D.R.
\emph{Special Functions}, World Scientific (1989).

\bibitem{krantz} Krantz, S.G.
\emph{Handbook of Complex Variables}, Springer Science, New York (1999).

\bibitem{prodanov} Prodanov, D., 
\emph{Regularized Integral Representations of the Reciprocal Gamma Function.} Fractal Fract \textbf{2019}, 3, 1.

\bibitem{reyn0}Reynolds, R.; Stauffer, A., 
\emph{Definite Integral of Arctangent and Polylogarithmic Functions Expressed as a Series.} Mathematics \textbf{2019}, 7, 1099.

\bibitem{reyn1}Reynolds, R.; Stauffer, A. 
\emph{A Definite Integral Involving the Logarithmic Function in Terms of the Lerch Function.} Mathematics \textbf{2019}, 7, 1148.

\bibitem{reyn2}Reynolds, R.; Stauffer, A. 
\emph{Derivation of Logarithmic and Logarithmic Hyperbolic Tangent Integrals Expressed in Terms of Special Functions.} Mathematics \textbf{2020}, 8, 687.

\bibitem{reyn3}Reynolds, R.; Stauffer, A. 
\emph{A Method for Evaluating Definite Integrals inTerms of Special Functions with Examples} International Mathematical Forum, Vol. 15, \textbf{2020}, no. 5, 235- 244

\bibitem{bdh} Bierens de Haan, D.,
\emph{Nouvelles Tables d'int\'egrales d\'efinies}, Amsterdam, 1867

\bibitem{choi} Choi, Junesang. and Srivastava, H.M.,
\emph{A family of log-gamma integrals and associated results}, Journal of Mathematical Analysis and Applications, 303 (2005) 436-449, Science Direct, Elsevier.

\bibitem{atlas}Jan Myland, Keith B Oldham, Jerome Spanier,
\emph{An Atlas of Functions: With Equator, the Atlas Function Calculator}, Springer; 2nd ed. 2009 edition (Dec 2 2008)

\bibitem{watson}Whittaker, E. T. and Watson, G. N. 
\emph{A Course in Modern Analysis}, 4th ed. Cambridge, England: Cambridge University Press. 

\bibitem{brychkov}Yu. A. Brychkov, O.I. Marichev, N.V. Savischenko,
\emph{Handbook of Mellin Transforms}, CRC Press, Taylor \& Francis Group, Boca Raton, FL., (2019)

\bibitem{prud} Prudnikov, A.P., Brychkov, Yu. A., Marichev, O.I.
\emph{Integrals and Series, More Special Functions}, USSR Academy of Sciences, Vol. 1, Moscow (1990).

\bibitem{as} Abramowitz, M. and Stegun, I.A.(Eds),
\emph{Handbook of Mathematical Functions with Formulas, Graphs, and Mathematical Tables, 9th printing}, New York, Dover, (1982).

\bibitem{malm1} Malmsten, C.J.,
\emph{De integralibus quibusdam definitis seriebusque infinitis (Eng. trans.: On some definite integrals and series)}, J. Reine Angew. Math., 38, 1-39 (1849) [work dated May 1, 1846].

\end{thebibliography}
\end{document}